\documentclass[a4paper,11pt]{amsart}
\usepackage{amssymb}
\usepackage{mathrsfs}
\usepackage{amsmath,amssymb,amsthm,latexsym,amscd,mathrsfs}
\usepackage{indentfirst}
 \newcommand{\ROM}[1]{\mathrm{\uppercase\expandafter{\romannumeral#1}}}
  \theoremstyle{definition}
 \newtheorem{theorem}{Theorem}[section]

 \newtheorem{corollary}{Corollary}[section]
  
 \newtheorem{remark}{Remark}[section]
 

\title[Gradient maps of isoparametric polynomials and its application]{\textbf{Gradient map of isoparametric polynomial and its application to Ginzburg-Landau system}}
\author[J. Q. Ge]{Jianquan Ge}\address{School of Mathematical Sciences, Laboratory of Mathematics and Complex Systems, Beijing Normal
University, Beijing 100875}\email{gejianquan@hotmail.com}
\thanks {The project is partially supported by the NSFC (Grant No. 10531090 and No.
10701007).}
\author[Y. Q. Xie]{Yuquan Xie}\address{School of Mathematical Sciences, Laboratory of Mathematics and Complex Systems, Beijing Normal
University, Beijing 100875}\email{yuqxie@gmail.com}
\thanks{The second author is the corresponding author.}
\date{}
\begin{document}
\maketitle
\begin{abstract}
In this note, we study properties of the gradient map of the
isoparametric polynomial. For a given isoparametric hypersurface in
sphere, we calculate explicitly the gradient map of its
isoparametric polynomial which turns out many interesting
phenomenons and applications. We find that it should map not only
the focal submanifolds to focal submanifolds, isoparametric
hypersurfaces to isoparametric hypersurfaces, but also map
isoparametric hypersurfaces to focal submanifolds. In particular, it
turns out to be a homogeneous polynomial automorphism on certain
isoparametric hypersurface. As an immediate consequence, we get the
Brouwer degree of the gradient map which was firstly obtained by
Peng and Tang with moving frame method. Following Farina's
construction, another immediate consequence is a counter example of
the Br\'ezis question about the symmetry for the Ginzburg-Landau
system in dimension 6, which gives a partial answer toward the Open
problem 2 raised by Farina.
\end{abstract}

%------------------------------------------------------------------------
\section{Introduction}
Let $M$ be a connected oriented isoparametric hypersurface in the
unit sphere $S^{n+1}$ with $g$ distinct principle curvatures. The
isoparametric polynomial $F$ of $M$ is a homogeneous polynomial of
degree $g$ in the Euclidean space $\mathbb{R}^{n+2}$, which is
uniquely determined by $M$ and satisfies the Cartan-M\"{u}nzner
equations:
\begin{eqnarray}
  &&|\nabla F|^2 =g^2|x|^{2g-2}, \label{eq1.1}\\
  &&\Delta F =\frac{g^2}{2}(m_2-m_1)|x|^{g-2}\label{eq1.2},
\end{eqnarray}
where $\nabla F$, $\Delta F$ denote the gradient and Laplacian of
$F$ in $\mathbb{R}^{n+2}$ respectively, and $m_1$, $m_2$ the
multiplicities of the maximal and minimal principal curvature of
$M$.

Cartan (see \cite{Ca38}, \cite{Ca39}) considered isoparametric
hypersurfaces in spheres and solved the classification problem in
the case $g\in \{1,2,3\}$. By using delicate cohomological and
algebraic arguments, M\"{u}nzner (see \cite{Mu80}, \cite{Mu81})
obtained the splendid result that the number $g$ must be $1,2,3,4$
or $6$. Deep going study about the geometry and topology of
isoparametric hypersurfaces leads to a lot of important results (see
\cite{TT72}, \cite{Tak76}, \cite{OT76}, \cite{Ab83}, \cite{DN85},
 \cite{Ta91}, \cite{Miy93}, \cite{CCJ07}, \cite{Imm08}, \emph{etc.}).

It is well known that there is a one-to-one correspondence between
isoparametric polynomials and families of parallel isoparametric
hypersurfaces. Throughout of this paper, we identify both
isoparametric polynomials and isoparametric hypersurfaces with their
congruences under isometries of the Euclidean space. Thus we say an
isoparametric polynomial is unique means that it's unique under
congruence. Given an isoparametric hypersurface $M$, one can
construct a function $F$ which turns out to be an isoparametric
polynomial with its level hypersurfaces being parallel hypersurfaces
of $M$. Given an isoparametric polynomial $F$ on $\mathbb{R}^{n+2}$,
let $f$ denote the restriction to $S^{n+1}$. Then the level
hypersurfaces of $f$ consist of a family of parallel isoparametric
hypersurfaces in the sphere. From the Cartan-M\"{u}nzner equation
$(\ref{eq1.1})$, it is not difficult to show that $f$ must have
$[-1,1]$ as its range. Moreover, the gradient of $f$ on $S^{n+1}$
can vanish only when $f=\pm 1$, and for each $s \in (-1,1)$, the
level set
\[   M_s=\{   p\in S^{n+1}| f(p)=s  \}   \]
is a compact connected isoparametric hypersurface, while $M_1$,
$M_{-1}$ are the focal submanifolds in the sphere with codimension
$m_1+1$ and $m_2+1$ respectively. In other words, the level sets of
$f$ give a ``singular" foliation of $S^{n+1}$ as
$S^{n+1}=\bigcup_{s\in [-1,1]}M_s $.

The gradient map $\Phi$ is a map from $\mathbb{R}^{n+2}$ to
$\mathbb{R}^{n+2}$ defined by $\Phi=\frac{1}{g}\nabla F$. Obviously,
each component of $\Phi$ is a homogeneous polynomial of degree
$g-1$, and the restriction of $\Phi$ to $S^{n+1}$ provides a
homogeneous polynomial map from $S^{n+1}$ to $S^{n+1}$. In
\cite{PT96}, Peng and Tang applied the moving frame method
successfully to obtain the Brouwer degree of $\Phi$. In this paper,
we try to study further properties of the gradient map of the
isoparametric polynomial and establish

\begin{theorem}
\label{th1.1}
  Let F be an isoparametric polynomial of degree $g$ on $\mathbb{R}^{n+2}$ ($g \geqslant 2$), the gradient map $\Phi=\frac{1}{g}\nabla F$.
Then $\Phi(M_{\cos(g\tau)})=M_{\cos(g(1-g)\tau)}$. In particular,
\begin{itemize}
  \item[(i)] $\Phi$ maps focal submanifold to focal submanifold;
  \item[(ii)] if $s \in \mathcal{D}=\{\cos \frac{k\pi}{g-1}| 1 \leqslant k < g-1\}$, $\Phi$ maps isoparametric hypersurface $M_s$ to
  focal submanifold;
  \item[(iii)] if $s \in (-1,1)-\mathcal{D}$, $\Phi$ maps isoparametric hypersurface $M_s$ to
  isoparametric hypersurface. In particular, $\Phi$ provides a homogeneous polynomial automorphism on certain isoparametric hypersurface.
\end{itemize}
\end{theorem}

As an immediate consequence, we can get the Brouwer degree of the
gradient map.

In Section 3, we give an application of the gradient maps to the
Open problem 2 of Farina \cite{Fa04} about Br\'{e}zis question on
the symmetry for the Ginzburg-Landau system.

\begin{bf}
Question (Br\'{e}zis \cite{Br99})
\end{bf} \quad Let $u: \mathbb{R}^N \rightarrow \mathbb{R}^N$ be a solution of
\begin{equation}\label{eq1.03}
  \Delta u=u(|u|^2-1)  \hspace{10pt} \texttt{on} \hspace{5pt}
\mathbb{R}^N ,\hspace{10pt} N\geqslant 3 \end{equation}
 with
$|u(x)|\rightarrow 1$ as $|x|\rightarrow \infty$ (possibly with a
``good" rate of convergence). Assume $\deg (u,\infty)=\pm 1$. Does
$u$ have the form
\begin{equation}
  u(x)=\frac{x}{|x|} h(|x|) \label{eq1.3}
\end{equation}
(modulo translation and isometry), where $h: \mathbb{R}_+\rightarrow
\mathbb{R}_+$ is a smooth function, such that $h(0)=0$ and
$h(\infty)=1$?

In \cite{Br99}, Br\'{e}zis gave an affirmative answer for the case
$N=2$ and thus raised the above question. For the case $N=8$, Farina
\cite{Fa04} gave a negative answer to it. In fact, he constructed a
radial solution which can be written as the form
$u(x)=G(\frac{x}{|x|})h(|x|)$. Therefore, he formulated his Open
problem 2 which is to study Br\'{e}zis question in dimension
$N\geqslant 3$ and $N\neq8$. Following Farina's construction and
using the gradient map of some isoparametric polynomial, we give
another counter example.
\begin{theorem}
\label{th1.2}
  There exists a solution $u: \mathbb{R}^6 \rightarrow
  \mathbb{R}^6$,
of the Ginzburg-Landau system (\ref{eq1.03}), satisfying
\begin{itemize}
  \item[(i)] $|u(x)|\rightarrow 1$ as $|x|\rightarrow \infty$;
  \item[(ii)] $\deg (u,\infty)=1$,
\end{itemize}
which has not the form $(\ref{eq1.3})$ (modulo translation and
isometry).

Furthermore, $u$ is a radial solution, \emph{i.e.}, it can be
written in the following way:
\begin{equation}
  u(x)=\Phi\left(\frac{x}{|x|}\right) h(|x|), \label{eq1.4}
\end{equation}
where $\Phi$ is the gradient map of the isoparametric polynomial in
$\mathbb{R}^6$ with $g=4, m_1=m_2=1$, and $h\in
C^2(\mathbb{R}_+,\mathbb{R})$ is the unique solution of
\begin{equation}
  \left\{
        \begin{array}{l}
          -h''-5\frac{h'}{r}+ 21\frac{h}{r^2}= h(1-h^2),~~ r>0,\\
          h(0)=0,\hspace{20pt} h(\infty)=1.
        \end{array}
  \right.
\end{equation}
\end{theorem}

\begin{remark}\label{rem1}
  Takagi \cite{Tak76} proved that for the isoparametric polynomial with $g=4$, if one of the
  principal curvatures of $M$ has multiplicity one, then $M$ must
  be homogeneous. Hence, the isoparametric polynomial in Theorem \ref{th1.2} is unique and one can write
  it as follows,
  \[   F=(|x|^2+|y|^2)^2-2\{(|x|^2-|y|^2)^2+ 4 \left<x,y\right>^2\},    \]
  where $(x,y)\in \mathbb{R}^3\times \mathbb{R}^3$ and $\left< , \right>$ denotes the
  Euclidean inner product in $\mathbb{R}^3$. Thus the map $\Phi=\frac{\nabla F}{4}$ in the theorem can be written explicitly.
\end{remark}

\begin{remark}\label{rem2}
  In the case $g=6, m_1=m_2=1$, Dorfmeister and Neher \cite{DN85}
  showed that the isoparametric hypersurface in $S^7$ must be homogeneous. In \cite{Miy93}, Miyaoka
  gave an interesting description in this case. She found that a homogeneous
  hypersurface in $S^7$ with $g=6$ is the inverse image of an
  isoparametric hypersurface in $S^4$ with $g=3$ under Hopf fibering.
  On the other hand, Cartan  \cite{Ca38}  \cite{Ca39} determined all isoparametric
  hypersurfaces with $g=3$. In particular, when $m_1=m_2=1$, the
  isoparametric polynomial in $ \mathbb{R}^5$ can be written as
  \[  F(x,y,X,Y,Z)= x^3-3xy^2+ \frac{3}{2}x(X^2+Y^2-2Z^2)+\frac{3\sqrt{3}}{2}y(X^2-Y^2)+3\sqrt{3}XYZ,  \]
  where $(x,y,X,Y,Z)\in\mathbb{R}^5$.
  Therefore, we can write the isoparametric polynomial with $g=6, m_1=m_2=1$
  as below
  \[\tilde{F}=F\circ \pi,\]
  where $\pi: \mathbb{R}^8 \rightarrow \mathbb{R}^5$ is given by
  \[   \pi(u,v)=(|u|^2-|v|^2,2u\bar{v}), \hspace{10pt} u,v \in \mathrm{the~quaternion~field} ~H\cong\mathbb{R}^4.    \]
  The gradient map $\Phi=\frac{\nabla \tilde{F}}{6}: \mathbb{R}^8\rightarrow\mathbb{R}^8$ is exactly the
  map $G$ in Farina's counter example as we mentioned before.
\end{remark}

\section{Gradient map of isoparametric polynomial}

In this section, following M\"unzner\cite{Mu80}, firstly we'll
construct the isoparametric polynomial from a given isoparametric
hypersurface (See also \cite{CR85}).

Suppose $M$ is a connected oriented isoparametric hypersurface in
$S^{n+1}$ with $g$ distinct principal curvatures
$\lambda_i:=\cot(\theta_i)$, $\lambda_1>\cdots>\lambda_g$. It's well
known that
\[  \theta_i=\theta_1+\frac{i-1}{g}\pi, \hspace{10pt} i=1,\cdots,g,  \]
and the multiplicity $m_i$ of $\lambda_i$ satisfy: $m_i=m_{i+2},~~
m_2=m_g$.

Let $\xi$ be the oriented unit normal vector field of $M$. Consider
the normal exponential map $\phi: M \times \mathbb{R} \rightarrow
S^{n+1}$ defined by
\begin{equation}
  \phi(x,t)= \cos t \cdot x+ \sin t \cdot \xi_x. \label{eq1.5}
\end{equation}
We know that $\phi$ has rank $n+1$, except where $\cot t$ is a
principal curvature of $M$. For any regular point (x,t) of $\phi$,
there exists an open neighborhood $U$ of $(x,t)$, such that $\phi$
is a diffeomorphism of $U$ onto $V=\phi(U)$. Define $\tau :
V\rightarrow \mathbb{R}$ by
\begin{equation}
  \tau (p)=\theta_1-t,\label{eq1.6}
\end{equation}
where $t$ is considered as a function on $V$ under the map $\phi$.
Note that the function $\tau$ is invariant if we take a parallel
isoparametric hypersurface of $M$ instead of $M$ in the definition.
Then we can obtain a homogeneous function $F$ on the cone of
$\mathbb{R}^{n+2}$ over $V$ by
\[  F(rp)=r^g \cos g\tau (p),\hspace{10pt} p\in V,~r>0.    \]
It is well known that the function $F$ is the restriction of a
homogeneous polynomial with degree $g$ in $\mathbb{R}^{n+2}$, which
is uniquely determined by $M$ and satisfies the Cartan-M\"{u}nzner
equations (\ref{eq1.1}) and (\ref{eq1.2}) (see \cite{Mu80},
\cite{CR85}). This polynomial is called the isoparametric polynomial
of $M$.

Let $\Phi$ be the (normalized) gradient map of $F$,  \emph{i.e.}
$\Phi=\frac{\nabla F}{g}$. Since for $g=1$, $\Phi=\nabla F$ is a
constant map. Hence in the following discussion, we assume that
$g=2,3,4$ or $6$.

Let $f$ denote the restriction of $F$ on $S^{n+1}$. For convenience,
we write the level set of $f$ as $\tilde{M}_{\tau}=f^{-1}(\cos g
\tau)$, $\tau\in \mathbb{R} $. It is not difficult to check that,
  \begin{enumerate}
    \item For any integer $j$,  $\tilde{M}_0=\tilde{M}_{\frac{2j\pi}{g}}$,
   $\tilde{M}_{\frac{\pi}{g}}=\tilde{M}_{\frac{(2j+1)\pi}{g}}$,  and $S^{n+1}=\bigcup_{\tau\in I_j}\tilde{M}_{\tau}$, where $I_j=[\frac{j\pi}{g},\frac{(j+1)\pi}{g}]$;
    \item $\tilde{M}_{0}$, $\tilde{M}_{\frac{\pi}{g}}$ are the focal submanifolds
    in $S^{n+1}$
    with codimension $m_1+1$, $m_2+1$, respectively;
    \item  For any $\tau \in (0,\frac{\pi}{g})$, $\tilde{M}_{\tau}$ is an isoparametric hypersurface with the maximal
          principal curvature $\cot \tau$,  \emph{i.e.} $\theta_1$ of
          $\tilde{M}_{\tau}$ equals $\tau$.
  \end{enumerate}

\begin{proof}[Proof of Theorem \ref{th1.1}.]
For convenience, let $M=f^{-1}(0)=\tilde{M}_{\frac{\pi}{2g}}$, then
$\theta_1=\frac{\pi}{2g}$ and
$\tilde{M}_{\tau}=\phi(M\times\{\frac{\pi}{2g}-\tau\})$. It is
easily seen that
\begin{equation}\label{normalvect}
\tilde{\xi}=-\sin (\frac{\pi}{2g}-\tau) \cdot x+ \cos
(\frac{\pi}{2g}-\tau)\cdot \xi
\end{equation}
is the oriented unit normal vector field of $\tilde{M}_{\tau}$.

We now calculate the gradient map $\Phi$. At each point
$p=\phi(x,\frac{\pi}{2g}-\tau) \in \tilde{M}_{\tau}\subset S^{n+1}$,
$\tau \in (0, \frac{\pi}{g})$,
\[   \nabla F(p)=\nabla_{_S} f(p) + \left<p,\nabla F\right> p,   \]
where $\nabla $, $\nabla_{S} $ are the gradient operators in
$\mathbb{R}^{n+2}$ and $S^{n+1}$ respectively.

Since $F$ is a homogeneous polynomial of degree $g$ and $f=\cos
g\tau $, by Euler's theorem,
\begin{equation}
  \nabla F(p)= -g \sin g\tau \cdot \nabla_{_S} \tau(p) + g F\cdot p.    \label{eq1.7}
\end{equation}
On the other hand, by (\ref{eq1.5}) and (\ref{eq1.6}), we have
\begin{equation}
  \nabla_{_S} \tau(p) =-\tilde{\xi},    \label{eq1.8}
\end{equation}
Substituting (\ref{eq1.8}) to (\ref{eq1.7}) implies that for each
$p=\phi(x,\frac{\pi}{2g}-\tau)$,
\begin{eqnarray}
  \Phi(p)&=&\frac{1}{g}\nabla F(p)= \cos g\tau(p) \cdot p+ \sin g\tau(p) \cdot \tilde{\xi}_p.    \label{eq1.10}\\
         &=& \cos(\frac{\pi}{2g}+(g-1)\tau)\cdot x + \sin(\frac{\pi}{2g}+(g-1)\tau)\cdot \xi_x   \nonumber\\
         &=& \phi(x,\frac{\pi}{2g}+(g-1)\tau), \nonumber
\end{eqnarray}
which follows that
\begin{eqnarray}
  \Phi(\tilde{M}_{\tau})= \tilde{M}_{(1-g)\tau}. \label{eq1.11}
\end{eqnarray}

By the continuity of the gradient map $\Phi$, equalities
$(\ref{eq1.10})$ and $(\ref{eq1.11})$ holds for all $\tau \in
[0,\frac{\pi}{g}]$. In particular, the focal submanifold
$\tilde{M}_0$ is the fix point set of $\Phi$. If $g$ is odd, $\Phi$
maps the other focal submanifold $\tilde{M}_{\frac{\pi}{g}}$ to
$\tilde{M}_{\frac{(1-g)\pi}{g}}=\tilde{M}_0$, and if $g$ is even,
$\Phi$ maps $\tilde{M}_{\frac{\pi}{g}}$ to
$\tilde{M}_{\frac{(1-g)\pi}{g}}=\tilde{M}_{\frac{\pi}{g}}$. Note
that $\Phi$ is just the antipodal map when restricted to the focal
submanifold $\tilde{M}_{\frac{\pi}{g}}$.

For $\tau\in (0,\frac{\pi}{g})$, it follows from the equality
$(\ref{eq1.11})$, $\Phi(\tilde{M}_{\tau})$ is a focal submanifold if
and only if $\cos (g(1-g)\tau)=\pm 1$, \emph{i.e},
\begin{equation*}
  \tau=\frac{k}{g(g-1)}\pi, \hspace{10pt} 1\leqslant k < g-1.
\end{equation*}

When $\cos(g\tau)=\cos(g(1-g)\tau)$, \emph{i.e}, $g=2$, or
$\tau=\frac{2k\pi}{g^2}$, or $\tau=\frac{2k\pi}{g^2-2g}$,
$k\in\mathbb{Z}$, $\Phi$ maps $\tilde{M}_{\tau}$ to itself. In
particular, there's always an isoparametric hypersurface with
$g\geqslant2$ on which $\Phi$ provides a homogeneous polynomial
automorphism. These complete the proof of Theorem $\ref{th1.1}$.
\end{proof}

As a consequence of Theorem $\ref{th1.1}$, we can deduce the Brouwer
degree of the gradient map $\Phi=\frac{\nabla F}{g}|_{S^{n+1}}:
S^{n+1}\rightarrow S^{n+1}$, by counting the number of inverse
points, counted with multiplicity $\pm1$ which is the sign of the
tangential map of $\Phi$ according whether it preserves or reverses
the orientation, of any regular value point. Our method here differs
from that of \cite{PT96} where they used the integral definition of
Brouwer degree and calculated it by moving frame method (See
 \cite{BT82} for different equivalent definitions of Brouwer
degree).

\begin{corollary}\label{cor2.1}
  Let $\Phi$ be the gradient map of an isoparametric polynomial with degree
  $g$. Then the Brouwer degree of $\Phi$ is given by
  \begin{enumerate}
    \item for g=2, $\deg \Phi= (-1)^{m_1+1}$;
    \item for g=3, $\deg \Phi= (-1)^{m_1+1}+(-1)^{m_1+m_2+1}$;
    \item for g=4, $\deg \Phi= (-1)^{m_1+1}+(-1)^{m_1+m_2+1}+(-1)^{m_2+1}$;
    \item for g=6, $\deg \Phi= 2\cdot
    (-1)^{m_1+1}+(-1)^{m_1+m_2+1}+(-1)^{m_2+1}-1$.
  \end{enumerate}
\end{corollary}

\begin{proof}
Denote by $J=(0,\frac{\pi}{g})$,
$J_k=(\frac{k-1}{g(g-1)}\pi,\frac{k}{g(g-1)}\pi)$, $1\leqslant
k\leqslant g-1$. Let $$\mathscr{M}:= S^{n+1}-(\tilde{M}_0\cup
\tilde{M}_{\frac{\pi}{g}})=\bigcup_{\tau\in J}\tilde{M}_{\tau},\quad
\mathscr{M}_k:=\bigcup_{\tau\in J_k}\tilde{M}_{\tau}.$$ Then
$\mathscr{M}$, $\mathscr{M}_k$ are open subsets of $S^{n+1}$ and
Theorem $\ref{th1.1}$ implies that $\Phi|_{\mathscr{M}_k}:
\mathscr{M}_k\rightarrow \mathscr{M}$ is a diffeomorphism for each
$1\leqslant k\leqslant g-1$. Thus every point $p$ in $\mathscr{M}$
is a regular value point of $\Phi$ and its inverse set equals
$\{p_k\in\mathscr{M}_k| \Phi(p_k)=p, k=1,\cdots,g-1.\}$ having $g-1$
points. Therefore, to calculate the Brouwer degree of $\Phi$, we
need only specify the sign of its tangential map $\Phi_*$ at each
$p_k$.

 Assume $p_k\in\tilde{M}_{\tau_k}$ for some $\tau_k\in J_k$.
 Then the principal curvatures of $\tilde{M}_{\tau_k}$ are given by
  $\lambda_{ki}=\cot(\tau_k +\frac{i-1}{g}\pi)$ with multiplicities $m_i$
  satisfying $m_i=m_{i+2}$, $m_2=m_g$, $i=1,\cdots,g$. Note that
  $m_1$, $m_2$ are determined by the isoparametric polynomial and thus
  are same for each $k$. Suppose $X$ is a principal tangent vector
  of $\tilde{M}_{\tau_k}$ with respect to $\lambda_{ki}$ at $p_k$.
  It is easily seen from formula $(\ref{eq1.10})$ that
\begin{equation}\label{tangent-t}
  \Phi_*(X)=\frac{\sin ((1-g)\tau_k +\frac{i-1}{g}\pi)}{\sin (\tau_k+\frac{i-1}{g}\pi) } X,
\end{equation}
where $X$ in the right side is regarded as the vector at $p$ by
parallel translating $X$ from $p_k$ to $p$ in $\mathbb{R}^{n+2}$. On
the other hand, from formulas $(\ref{normalvect})$, $(\ref{eq1.8})$
and $(\ref{eq1.10})$, we can derive directly
\begin{equation}\label{tangent-n}
 \Phi_*(\tilde{\xi}_{p_k})=(1-g)\tilde{\xi}_p,
\end{equation}
where $\tilde{\xi}_{p_k}$, $\tilde{\xi}_p$ are the unit normal
vectors of $\tilde{M}_{\tau_k}$,
$\Phi(\tilde{M}_{\tau_k})=\tilde{M}_{(1-g)\tau_k}$ at $p_k$ and $p$
respectively. Notice that the tangent space of $S^{n+1}$ at $p_k$
(\emph{resp.} $p$) is spanned by such $X$s and $\tilde{\xi}_{p_k}$
(\emph{resp.} $\tilde{\xi}_p$), it follows immediately from
$(\ref{tangent-t})$ and $(\ref{tangent-n})$ that the sign of the
tangential map $\Phi_*$ at $p_k$ is given by
\begin{equation}
  \mathrm{sign~}
  \Phi_*|_{p_k}=
  \begin{cases}
    (-1)^{\frac{k+1}{2}m_1+\frac{k-1}{2}m_2+1} & \mathrm{for~} $k$
    \mathrm{~is~odd,} \\
    (-1)^{\frac{k}{2}m_1+\frac{k}{2}m_2+1} & \mathrm{for~} $k$
    \mathrm{~is~even.}
  \end{cases}\label{eq1.17}
\end{equation}

Combining $(\ref{eq1.17})$ with the following formula for the
definition of Brouwer degree
\[  \deg (\Phi)=\sum\limits_{k=1}^{g-1} \mathrm{sign~}
  \Phi_*|_{p_k},
\]
we can conclude the items of Corollary \ref{cor2.1}.
\end{proof}

When the isoparametric polynomial $F$ is harmonic, \emph{i.e.}
$m_1=m_2=:m$. According to the tangential map $\Phi_*$ given in
$(\ref{tangent-t})$ $(\ref{tangent-n})$, one can calculate directly
that the tension field $B(\Phi):=\mathrm{Trace}(\nabla_S\Phi_*)=0$ ,
hence $\Phi|_{S^{n+1}}:S^{n+1}\rightarrow S^{n+1}$ is a harmonic map
(See also \cite{ER93} and \cite{PT96}). For applications in next
section, we now focus on harmonic isoparametric polynomials.

For $g=2$, the harmonic isoparametric polynomial is given by
$F(x,y)=|x|^2-|y|^2$, and thus $\Phi(x,y)=(x,-y)$, where
$(x,y)\in\mathbb{R}^{m+1}\times\mathbb{R}^{m+1}$.  For $g=3$, Cartan
completely classified the isoparametric polynomials and showed that
$m_1= m_2=1,2,4, or~ 8$. For g=4, Abresch \cite{Ab83} showed that
harmonic isoparametric polynomials must have $m_1=m_2=1, or~2$.
These two cases were showed to be unique by Ozeki and Takeuchi
\cite{OT76}. See Remark \ref{rem1} for explicit representation of
the one with $m=1$. For $g=6$, M\"{u}nzner (see \cite{Mu80},
\cite{Mu81}) showed that it must have $m_1=m_2$. Furthermore,
Abresch \cite{Ab83} was able to show that the common multiplicity
$m$ must be either 1 or 2. In the case $m=1$, Dorfmeister and Neher
showed that it must be homogeneous. See Remark \ref{rem2} for
explicit representation for the case of $m=1$. Recently, Miyaoka
\cite{Miy08} claimed that it is also unique for the case of $m=2$.

In conclusion, by Corollary \ref{cor2.1} and discussions above, we
have (Compare with \cite{Ta07})
\begin{corollary}
\label{cor2.2} Harmonic isoparametric polynomial exists only when
$$(g,m)=(1,m),(2,m),(3,1),(3,2),(3,4),(3,8),(4,1),(4,2),(6,1),(6,2),$$
and for each case (except possibly the last one), it's unique under
congruence. Furthermore, its gradient map
$\Phi|_{S^{n+1}}:S^{n+1}\rightarrow S^{n+1}$ is a polynomial
harmonic map with the Brouwer degree $\deg \Phi$
  \begin{enumerate}
    \item $(g,m)=(1,m)$,\quad $\deg \Phi= 0$;
    \item $(g,m)=(2,m)$,\quad $\deg \Phi= (-1)^{m+1}$;
    \item $(g,m)=(3,1)$, \quad$\deg \Phi= 0$,\quad \quad$(g,m)=(3,2),(3,4),(3,8)$, \quad$\deg \Phi= -2$;
    \item $(g,m)=(4,1)$, \quad$\deg \Phi= 1$,\quad \quad$(g,m)=(4,2)$, \quad $\deg \Phi=-3$;
    \item $(g,m)=(6,1)$, \quad$\deg \Phi= 1$,\quad \quad$(g,m)=(6,2)$, \quad $\deg \Phi=-5$.
  \end{enumerate}
\end{corollary}

\section{Proof of Theorem \ref{th1.2}}

The proof should be translated word by word from Farina \cite{Fa04}
once one knows examples of harmonic isoparametric polynomial with
Brouwer degree of its gradient map being $\pm1$. For completeness,
we state it as follows.

By the results of \cite{AF97} a function $u$ having the form
$(\ref{eq1.4})$: $u(x)=\Phi(\frac{x}{|x|})h(|x|)$, with a
non-constant $\Phi \in C^2(S^{N-1}, \mathbb{R}^N)$ and a profile $h
\in C^2(\mathbb{R}_+,\mathbb{R})$ is a solution of the
Ginzburg-Landau system $(\ref{eq1.03})$ if
  \begin{enumerate}
    \item  $\Phi(S^{N-1})\subset S^{N-1}$,
    \item there exists a positive integer $k$ such that $\Phi \in
    (\mathscr{S}\mathscr{H}_{k,N})^N$ (where
    $\mathscr{S}\mathscr{H}_{k,N}$ is the vector space of the spherical
    harmonics of degree $k$ in $\mathbb{R}^N$ ),
  \end{enumerate}
and the profile $h$ satisfies
\begin{equation}
  \left\{
        \begin{array}{l}
          -h''-(N-1)\frac{h'}{r}+ k(k+N-2)\frac{h}{r^2}= h(1-h^2), r>0,\\
          h(0)=0.
        \end{array}
  \right. \label{eq1.14}
\end{equation}

Therefore, to obtain the desired conclusion it is enough to prove
the existence of a map $\Phi$ satisfying (i) and (ii) above with
$N=6$ and $k=3$, and a corresponding profile $h$ satisfying
$(\ref{eq1.14})$ with $h(\infty)=1$.

\emph{Existence of} $\Phi$: Corollary $\ref{cor2.2}$ implies that
the gradient map $\Phi=\frac{\nabla F}{g}$ of a harmonic
isoparametric polynomial $F$ has Brouwer degree $\pm1$ if and only
if $g=2$, or $(g,m)=(4,1)$, or $(g,m)=(6,1)$. Obviously, such map
$\Phi$ satisfies properties (i) and (ii) above. As mentioned before,
when $g=2$, $\Phi(x,y)=(x,-y)$ is congruent to the identity and so
is trivial. The map $\Phi$ for $(g,m)=(6,1)$ is given explicitly in
Remark \ref{rem2} and has been applied in Farina's counter example
in dimension $N=gm+2=8$. The map $\Phi$ for $(g,m)=(4,1)$ is exactly
the one we apply to construct the counter example in dimension
$N=6$. See Remark \ref{rem1} for explicit form of this map $\Phi$.

\emph{Existence of} $h$: In  \cite{FG00} it is proved that there is
a unique solution of the problem
\begin{equation}
  \left\{
        \begin{array}{l}
          -h''-(N-1)\frac{h'}{r}+ k(k+N-2)\frac{h}{r^2}= h(1-h^2), r>0,\\
          h(0)=0,\hspace{20pt} h(\infty)=1
        \end{array}
  \right. \label{eq1.15}
\end{equation}
for every integer $N \geqslant 3$ and every positive integer $k$.
Furthermore, the profile $h$ is a strictly increasing function. This
property implies that $u$ satisfies the condition (i) in Theorem
$\ref{th1.2}$. On the other hand, owing to the special form of the
constructed radial solution $u$, we have that $\deg (u,\infty)$ is
equal to the Brouwer degree of the map $\Phi$. This shows that (ii)
in Theorem $\ref{th1.2}$ is also satisfied. This completes the proof
of Theorem $\ref{th1.2}$.
\bigskip\\
{\bf Acknowledgements} It's our great pleasure to thank Professors
Chiakuei Peng and Zizhou Tang for introducing this topic to us and
guidance, and also for their careful revision of an earlier version
of this paper.

\end{document}